\documentclass{MIPRO}

\usepackage{cite}
\usepackage{amsmath,amssymb,amsfonts}
\usepackage{algorithmic}
\usepackage{graphicx}
\usepackage{textcomp}
\usepackage{xcolor}
\usepackage[T1]{fontenc}
\usepackage{flushend}
\usepackage{multirow}
\usepackage{multicol}
\usepackage{array}
\usepackage{url}
\usepackage{xcolor}

\newcommand{\fdmdelta}{\delta}

\DeclareMathOperator{\id}{id}

\begin{document}

\title{A Numerical Study of Combining RBF Interpolation and Finite Differences to Approximate Differential Operators}

\author{
\IEEEauthorblockN{
Adrijan Rogan\IEEEauthorrefmark{1}\IEEEauthorrefmark{2},
Andrej Kolar-Požun\IEEEauthorrefmark{1}\IEEEauthorrefmark{2},
Gregor Kosec\IEEEauthorrefmark{1}
}

\IEEEauthorblockA{\IEEEauthorrefmark{1} 
Parallel and Distributed Systems Laboratory, Jožef Stefan Institute, Ljubljana, Slovenia}

\IEEEauthorblockA{\IEEEauthorrefmark{2} 
Faculty of Mathematics and Physics, University of Ljubljana, Ljubljana, Slovenia}

adrijan.rogan9@gmail.com
}

\maketitle

\begin{abstract}
This paper focuses on RBF-based meshless methods for approximating differential operators, one of the most popular being RBF-FD. Recently, a hybrid approach was introduced that combines RBF interpolation and traditional finite difference stencils. We compare the accuracy of this method and RBF-FD on a two-dimensional Poisson problem for standard five-point and nine-point stencils and different method parameters.
\end{abstract}

\renewcommand\IEEEkeywordsname{Keywords}

\section{Introduction}
The forthcoming study focuses on numerical approximation of differential operators. Efficient high order approximations can be obtained by the use of Finite Difference (FD) schemes, however, these operate on structured, grid-based node layouts, limiting their flexibility in the setting of complex geometries.

For that reason, approximation methods that can operate on irregularly positioned or scattered nodes are of considerable interest. A common approach for such problems are mesh-based methods, such as the finite element or the finite volume method. The focus of this paper, however, is an alternative approach -- meshless methods, based on Radial Basis Functions (RBFs)~\cite{buhmann2000}.

RBFs have first appeared in the 70s in the context of scattered data interpolation~\cite{rbf70s} and have since enjoyed increasing popularity due to their provable invertibility guarantees and dimensionality independence.

In the decades that followed they have found application in various fields, including computer graphics, machine learning, finance and numerical solutions of Partial Differential Equations (PDEs)~\cite{fasshauer-matlab2007}.

In our study we focus on applications of RBFs in the latter - we aim to approximate linear differential operators on scattered node layouts. RBFs have been first used in this context by Kansa in 1990~\cite{kansa}, resulting in a global collocation method for meshless solutions of PDEs. In the 2000s~\cite{rbffd} a local version of the method, known as the Radial Basis Function generated Finite Differences (RBF-FD) appeared, directly generalizing the well-known FD to a scattered setting -- much like the usual FD, RBF-FD approximation can be obtained by starting with an interpolant and applying the desired differential operator to it.

Another way to generalize FD to a meshless setting is to use RBFs to interpolate the function values to a "{}virtual"{} finite difference stencil and then apply the usual FD to evaluate the derivative. Such hybrid approach was recently explored and successfully applied on an elasto-plasticity problem, where it was demonstrated that it can outperform RBF-FD~\cite{vuga2024}. A potential benefit of this approach is the possibility of adapting existing FD schemes to a meshless setting.

The goal of this paper is to compare the accuracy of the two approaches on a model 2D Poisson problem. In the following Section we briefly describe RBF interpolation and its role in approximating differential operators. In Section 3 we then present our problem setup for the comparison of the two methods. In Section 4 our main results are presented and then summarised in Section 5.

\section{Methods}

\subsection{RBF interpolation}
\begin{figure}
  \centering
  \includegraphics[width=0.48\textwidth]{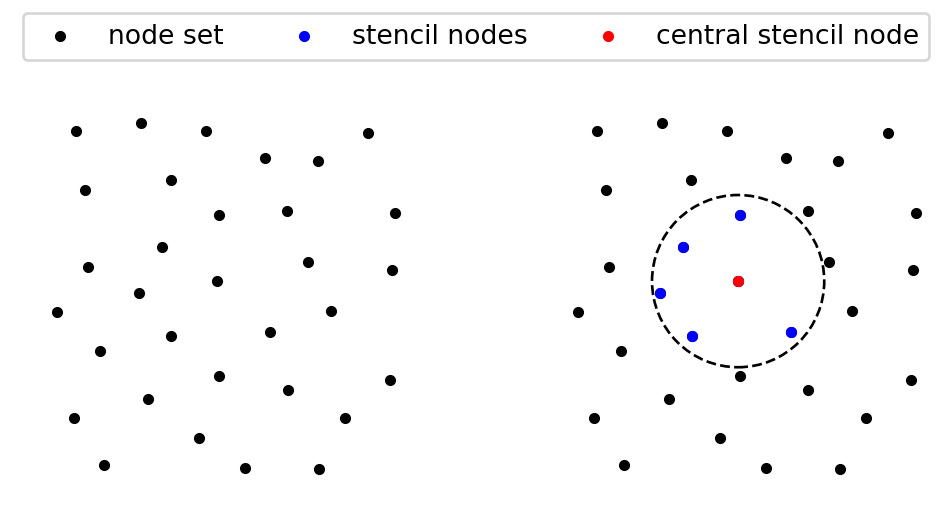}
  \caption{An example of a scattered node set $X \subset \mathbb{R}^2$ on the left-hand side and a stencil with six nearest nodes to the central node on the right-hand side. For clarity, stencil nodes are additionally encircled.}
  \label{fig:plot_domain}
\end{figure}

Consider a scattered node set $X = \{{\boldsymbol{x}_i}\}_{i=1}^N$ with corresponding function values $f_i = f(\boldsymbol{x_i})$ and a chosen function $\phi$. An example node set is presented on the left-hand side of Fig. \ref{fig:plot_domain}. The (global) RBF interpolant takes the form
\begin{equation}
    \hat{f}(\boldsymbol{x}) = \sum_{i=1}^N \alpha_i \phi(\Vert \boldsymbol{x} - \boldsymbol{x_i} \Vert).
\end{equation}

Defining $\varphi_i(\boldsymbol{x}) = \phi(\Vert \boldsymbol{x} - \boldsymbol{x_i} \Vert)$, we obtain a system of linear equations from the interpolation conditions:
\begin{equation}
    \label{eqn:interpolation-matrix}
    \begin{bmatrix}
         \varphi_1(\boldsymbol{x_1}) & \cdots & \varphi_N(\boldsymbol{x_1}) \\ \vdots & \ddots & \vdots \\ \varphi_1(\boldsymbol{x_N}) & \cdots & \varphi_N(\boldsymbol{x_N})
    \end{bmatrix}
    \begin{bmatrix}
        \alpha_1 \\ \vdots \\ \alpha_N
    \end{bmatrix}
    = 
    \begin{bmatrix}
        f_1 \\ \vdots \\ f_N
    \end{bmatrix}
    .
\end{equation}
This system, written compactly as $A\boldsymbol{\alpha} = \boldsymbol{f}$ can be proven to be uniquely solvable if $\phi$ is a positive definite RBF and the nodes are pairwise distinct. For non-compactly supported RBFs, the matrix $A$ is dense and becomes ill-conditioned for large node sets \cite{fasshauer-matlab2007}.

To address the problem of ill-conditioning, we can consider local interpolation around an arbitrary position $\boldsymbol{\hat{x}}$, that is not necessarily part of $X$. To do this, we introduce the stencil\footnote{Another common term is support.} $S(\boldsymbol{\hat{x}})$ as a set of indices corresponding to a subset of $X$, usually taken to be a fixed number $n$ of closest nodes to $\boldsymbol{\hat{x}}$. The term stencil often refers directly to the node subset and not to the index set. An example is shown on Fig. \ref{fig:plot_domain}. The local interpolant is then:
\begin{equation}
\label{eqn:interpolation-local}
    \hat{f}(\boldsymbol{x}) = \sum_{i \in S(\boldsymbol{\hat{x}})} \alpha_i \phi(\Vert \boldsymbol{x} - \boldsymbol{x_i} \Vert).
\end{equation}

This results in a smaller linear system $A\boldsymbol{\alpha} = \boldsymbol{f}$ (limited by the stencil size). The interpolant is often augmented with monomials:
\begin{equation}
    \label{eqn:rbf-augmentation}
    \hat{f}(\boldsymbol{x}) = \sum_{i \in S(\boldsymbol{x})} \alpha_i \phi(\Vert \boldsymbol{x} - \boldsymbol{x_i} \Vert) + \sum_{i=1}^s \beta_i p_i(\boldsymbol{x}),
\end{equation}
where $p_i$ are different monomials. It is customary to include all monomials up to some degree $m$ inclusive, which results in $s= \binom{m+2}{m}$ if $X \subset \mathbb{R}^2$. This ensures the interpolant has polynomial reproduction of the same degree, but with increased stability compared to the usual, purely polynomial approximation~\cite{bayona2017}. In this case, we must solve an extended linear system with $M = \vert S(\boldsymbol{x}) \vert + s$ equations and $M$ unknowns. For further details, see chapter 6 of \cite{fasshauer-matlab2007}.

\subsection{The RBF-FD Method}

Let $\boldsymbol{x}$ be a node at which we want to approximate a differential operator $\mathcal{L}$. We seek an approximation in the following form:
\begin{equation}
    \label{eqn:rbffd-operator-approximation}
    (\mathcal{L}u)(\boldsymbol{x}) \approx \sum_{i \in S(\boldsymbol{x})} w_i u(\boldsymbol{x_i}) = \boldsymbol{w}^\top \boldsymbol{u}.
\end{equation}

To determine the weights $\boldsymbol{w}$, we construct a local RBF interpolant $\hat{u}$. Following (\ref{eqn:interpolation-local}), the interpolant weights are $\boldsymbol{\alpha} = A^{-1}\boldsymbol{u}$. Applying the differential operator, we obtain:
\begin{align}
    (\mathcal{L}u)(\boldsymbol{x}) &\approx (\mathcal{L}\hat{u})(\boldsymbol{x}) = \sum_{i \in S(\boldsymbol{x})} \alpha_i (\mathcal{L}\varphi_i)(\boldsymbol{x}) = \\
    & = (\mathcal{L}\boldsymbol{\varphi})(\boldsymbol{x})^\top \boldsymbol{\alpha} = (\mathcal{L}\boldsymbol{\varphi})(\boldsymbol{x})^\top A^{-1} \boldsymbol{u} \nonumber,
\end{align}
therefore the weights $\boldsymbol{w}$ are the solution of the linear system $A \boldsymbol{w} = (\mathcal{L}\boldsymbol{\varphi})(\boldsymbol{x})$.

As in interpolation, monomial augmentation is often used when constructing the local interpolant. In this case, the form of the approximation (\ref{eqn:rbffd-operator-approximation}) stays the same but determining the weights $\boldsymbol{w}$ again requires solving an extended linear system. More information can be found in \cite{flyer2016}.

The computational complexity of obtaining the weights $\boldsymbol{w}$ depends on the method used to solve the linear system. We used an LU decomposition with partial pivoting, which results in an asymptotic flop count of $\frac{2}{3}M^3$. Additionally, computing the solution with the decomposition requires one forward and one backward substitution, each taking $M^2$ flops \cite{numlinalg1997}.

\subsection{Combining RBF interpolation and finite differences}

Finite difference schemes typically use a uniform grid for the domain discretization and differential operators are approximated with finite difference stencils. These stencils are operator-dependent and are given by their offsets $\{\Delta_1, \dots, \Delta_k\}$ and weights $\{a_1, \dots, a_k\}$. A given operator $\mathcal{L}$ is approximated at position $\boldsymbol{x}$ as: 
\begin{equation}
    \label{eqn:fdm-approx}
    (\mathcal{L}u)(\boldsymbol{x}) \approx \sum_{i=1}^k a_i u(\boldsymbol{x} + \Delta_i).
\end{equation}


We can try and apply finite difference stencils to a scattered node set; we will refer to them as virtual stencils in this context. The approximation of a differential operator is the same as in (\ref{eqn:fdm-approx}) but in general we cannot use it directly as some virtual stencils may refer to points outside the node set. We can circumvent this by interpolating
\begin{equation}
    \label{eqn:hybrid-interp}
    u(\boldsymbol{x}+\Delta_i) \approx \sum_{j \in S(\boldsymbol{x})} w_{ij} u(\boldsymbol{x_j})
\end{equation}
to obtain the approximation
\begin{equation}
    \label{eqn:hybrid-approx}
    (\mathcal{L}u)(\boldsymbol{x}) \approx \sum_{i=1}^k a_i u(\boldsymbol{x}+\Delta_i) \approx \sum_{i=1}^k a_i \sum_{j \in S(\boldsymbol{x})} w_{ij} u(\boldsymbol{x_j}),
\end{equation}

\begin{figure}
  \centering
  \includegraphics[width=0.48\textwidth]{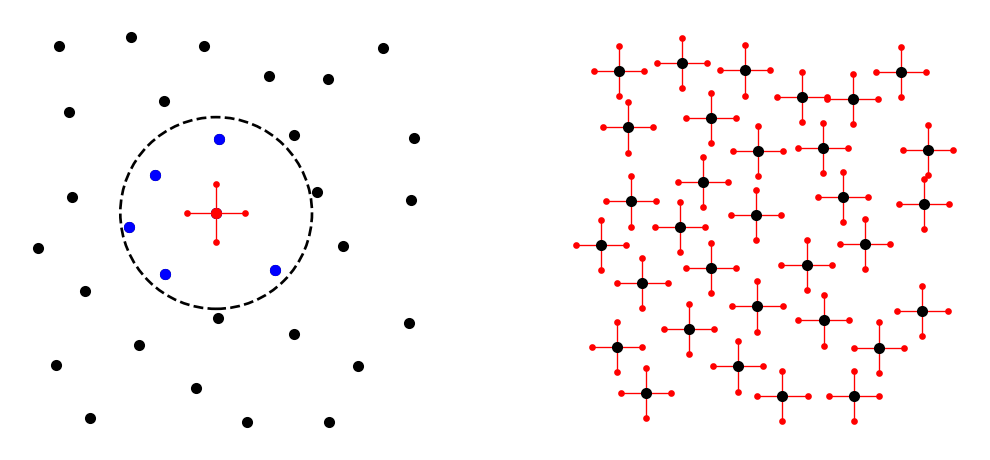}
  \caption{An example of a virtual stencil for a single point on the left-hand side and for all points from the domain on the right-hand side.}
  \label{fig:plot_stencil}
\end{figure}
An example of a scattered node set $X$ along with virtual stencils is shown on Fig. \ref{fig:plot_stencil}.
Note that (\ref{eqn:hybrid-interp}) can be viewed as a special case of RBF-FD for $\mathcal{L} = \id$ (the identity operator) and we determine the weights in a similar way. The weights depend only on the positions of the support nodes, therefore they can be computed only once and stored.

An important aspect to consider is computation time, which is higher compared to the RBF-FD method: The weights $\{a_1, \dots, a_k\}$ in (\ref{eqn:hybrid-approx}) are given, but, the weights $w_{ij}$ are determined by solving $k$ linear systems\footnote{If one of the stencil offsets $\Delta_i$ is zero, interpolation (\ref{eqn:hybrid-interp}) is unnecessary for that offset, so there are $k$ linear systems to solve} of size $M \times M$. However, the matrix for each stencil offset is the same, so this can be done efficiently if we store the LU decomposition of the matrix, leaving the leading term unchanged at $\frac{2}{3}M^3$ flops. Nevertheless, we still need to obtain the solutions of $k$ linear systems instead of just one. Therefore, this approach requires additional $2(k-1)M^2$ flops to perform $k$ forward and backward substitutions compared to RBF-FD.

\section{Problem setup}
We consider the Poisson problem on a unit square domain $\Omega = (0,1)^2$, i.e.
\begin{align}
    \label{eqn:poisson-problem}
    \Delta u &= f \; \text{in } \Omega \\ 
    u &= 0 \; \text{on } \partial \Omega, \nonumber
\end{align}
where the right-hand side function is given by
\begin{equation}
    f(x,y) = -2\pi^2\sin(\pi x)\sin(\pi y),
\end{equation}
which corresponds to the analytical solution
\begin{equation}
    u(x,y) = \sin(\pi x)\sin(\pi y).
\end{equation}
We solve this problem using both the RBF-FD method and the hybrid approach described earlier, and compare the results.

For both methods, we use the same set of discretization points $X$. The boundary of the domain is first discretized with a uniform step size $h$ and the interior is filled with fill density $h$ using the algorithm described in \cite{slak2019}. We limit our testing to a fixed value of $h= 0.01$ and use a fixed seed for the random number generator so that the node set $X$ remains unchanged throughout testing.

We opt for a popular choice of a polyharmonic spline RBF of order $k=3$, i.e. $\phi(r)=r^3$, which has been shown to have desirable properties regarding accuracy and stability under the condition that we augment with monomials of sufficiently high degree \cite{flyer2016,bayona2017}. We will test different orders of augmentation $m \in \{2, 4, 6, 8\}$.

Each point in the discretization has a corresponding support consisting of its $n$ nearest neighbors (including itself). Following~\cite{bayona2017}, we set $n$ to be twice the number of augmenting monomials:
\begin{equation}
    n = 2 \cdot \binom{m+2}{m}.
\end{equation}

For the hybrid method, we approximate the Laplace operator via approximation of second-order derivatives in each spatial direction. A commonly used finite difference stencil for the second-order derivative is
\begin{equation}
    u_{xx}(x) \approx \frac{1}{\fdmdelta^2} \bigl( u(x-\fdmdelta) - 2u(x) - u(x+\fdmdelta) \bigr),
\end{equation}
where $\fdmdelta$ is some chosen spacing of the FD stencil.
Applying this to both spatial dimensions, we approximate the Laplacian using a five-point stencil:
\begin{multline}
    \Delta u(x, y) = u_{xx}(x,y) + u_{yy}(x,y) \\ \approx \frac{1}{\fdmdelta^2} \cdot \bigl(u(x+\fdmdelta,y) + u(x,y+\fdmdelta) + u(x-\fdmdelta,y) \\ + u(x, y-\fdmdelta) - 4u(x,y) \bigr).
\end{multline}
The five-point stencil is exact for polynomials up to degree $3$, inclusive. For improved accuracy, we can also use higher order stencils. In addition to the five-point stencil, we have also tested a nine-point stencil, which is exact for polynomials up to degree $5$, inclusive.
As the polynomial reproduction degree of an RBF interpolant, augmented with monomials of degree $m$ is equal to $m$, the polynomial reproduction degree of the hybrid method is $\mathrm{min}(m,3)$ for the five-point stencil and $\mathrm{min}(m,5)$ for the nine-point stencil.
The weights and offsets of both stencils are presented in Table \ref{tab:stencil}. We will vary the virtual stencil spacing $\fdmdelta$ by sampling $\sigma \in \mathbb{R}^+$ and setting $\fdmdelta = \sigma \cdot h$.

\begin{table*}
    \centering
    \begin{tabular}{c|c|c|c|c|c|c|c|c|c}
         index $i$ & $1$ & $2$ & $3$ & $4$ & $5$ & $6$ & $7$ & $8$ & $9$ \\ \hline
         offset $\Delta_i$ & $(0, 0)$ & $(\fdmdelta, 0)$ & $(-\fdmdelta, 0)$ & $(0, \fdmdelta)$ & $(0, -\fdmdelta)$ & $(2\fdmdelta, 0)$ & $(-2\fdmdelta, 0)$ & $(0, 2\fdmdelta)$ & $(0, -2\fdmdelta)$ \\ \hline
         5-point stencil weight $\fdmdelta^2 \cdot a_i$ & $-4$ & $1$ & $1$ & $1$ & $1$ & / & / & / & / \\ \hline
         9-point stencil weight $\fdmdelta^2 \cdot a_i$ & $-5$ & $4/3$ & $4/3$ & $4/3$ & $4/3$ & $-1/12$ & $-1/12$ & $-1/12$ & $-1/12$ \\
    \end{tabular}
    \caption{Weights and offsets of the virtual five-point and nine-point stencils.\\The five-point stencils use only the first five offsets.}
    \label{tab:stencil}
\end{table*}

The solutions obtained from both RBF-FD and the hybrid approach are approximate values $\tilde{u}(\boldsymbol{x})$ at each discretization point $\boldsymbol{x} \in X$. We obtain $\tilde{u}(\boldsymbol{x_1}), \dots, \tilde{u}(\boldsymbol{x_N})$ by solving a sparse linear system, where the matrix is assembled with the weights from the operator approximation at each node (much like in the usual finite difference method) and the right-hand side is given by (\ref{eqn:poisson-problem}).

We compare their accuracy by computing the mean and the maximum relative errors with respect to the analytical solution $u(\boldsymbol{x})$, limited to interior nodes:
\begin{equation}
    \text{maximum relative error} = \max_{\boldsymbol{x} \in \Omega \cap X} \left| \frac{u(\boldsymbol{x}) - \tilde{u}(\boldsymbol{x})}{u(\boldsymbol{x})} \right|,
\end{equation}
\begin{equation}
    \text{mean relative error} = \frac{1}{N} \cdot \sum_{\boldsymbol{x} \in \Omega \cap X} \left| \frac{u(\boldsymbol{x}) - \tilde{u}(\boldsymbol{x})}{u(\boldsymbol{x})} \right|.
\end{equation}

Additionally, for both methods, we measure the runtime of two phases of the solution procedure: first, the construction of the large sparse linear system (i.e., weight computations), and second, solving the system.

All methods and solution procedures were implemented using the Medusa library \cite{medusa} and the code is publicly available\footnote{\url{https://gitlab.com/e62Lab/2025_cp_mipro_combiningrbfandfdm}}.

\section{Results}
\begin{figure}
    \centering
    \includegraphics[width=\linewidth]{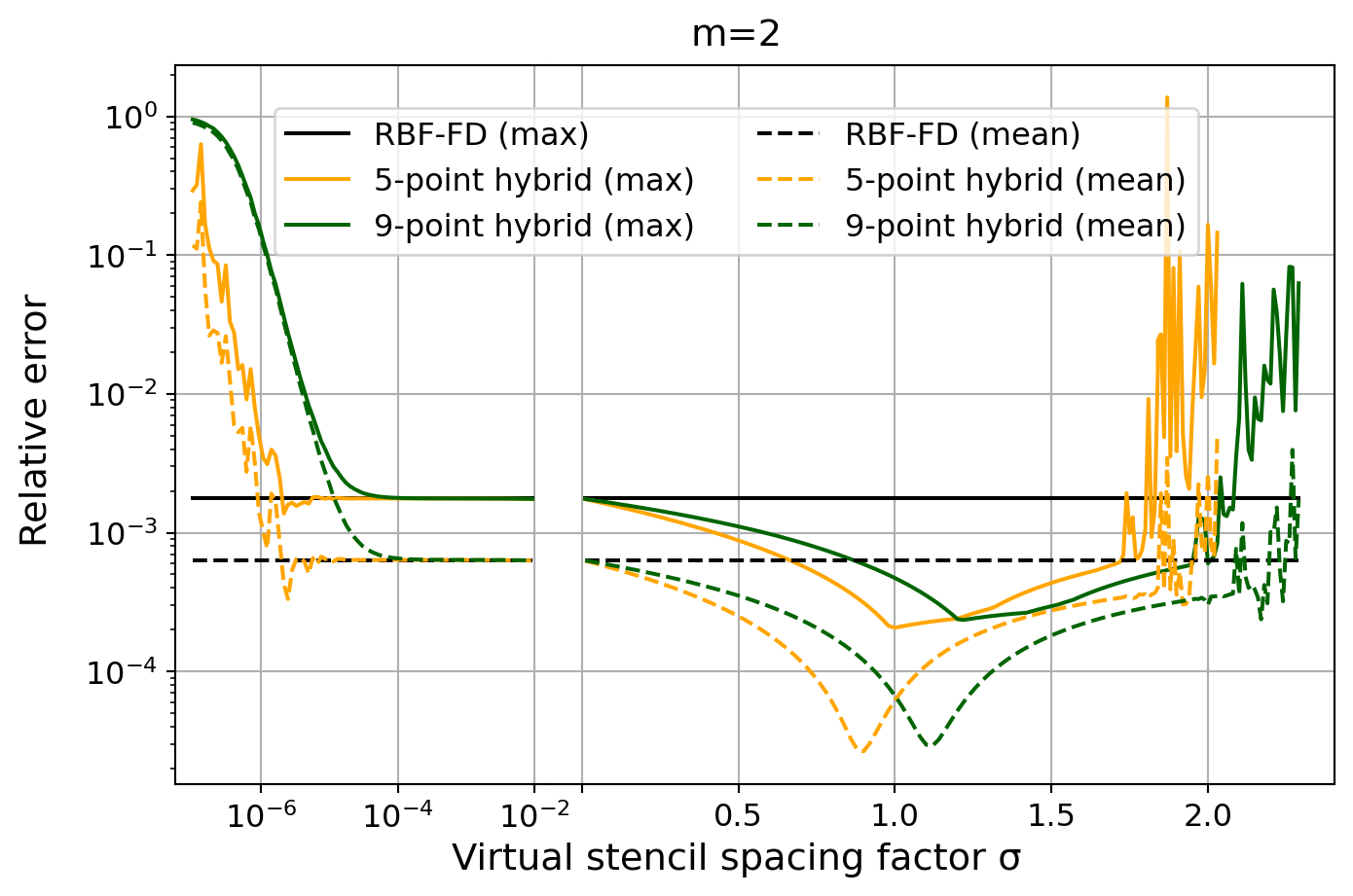}
    \caption{Error comparison for 2nd order monomial augmentation and the five-point and nine-point stencil.}
    \label{fig:error-m2}
\end{figure}

We begin our analysis with monomial augmentation of order $m=2$ and vary the virtual stencil size by varying $\sigma$. The results can be seen on Fig. \ref{fig:error-m2}. We notice that for small values of $\sigma$, the error of the hybrid approach is similar to the RBF-FD error, but starts to diverge as $\sigma$ becomes increasingly small. This is expected already because of the limited numerical precision - $\sigma^2$ appears in the denominator of the virtual stencil. Likewise, the error grows and becomes noisy also for large values of $\sigma$, likely because the virtual stencil nodes become too distant from the stencil center, affecting the accuracy of the interpolation step. In fact, for the larger values of $\sigma$ issues when solving the global sparse system started to arise. Between these two extremes, there is a sweet spot of the hybrid approach at approximately $\sigma \approx 1$, where the virtual stencil spacing $\delta$ is about the same as the fill density $h$. Interestingly, the hybrid approach at this point performs much better than RBF-FD, reducing the relative error by up to an order of magnitude. Note also that the nine-point stencil does not provide an improvement over the five-point stencil, which is expected as the error is dominated by the less accurate $m=2$ interpolation step.

\begin{figure}
    \centering
    \includegraphics[width=\linewidth]{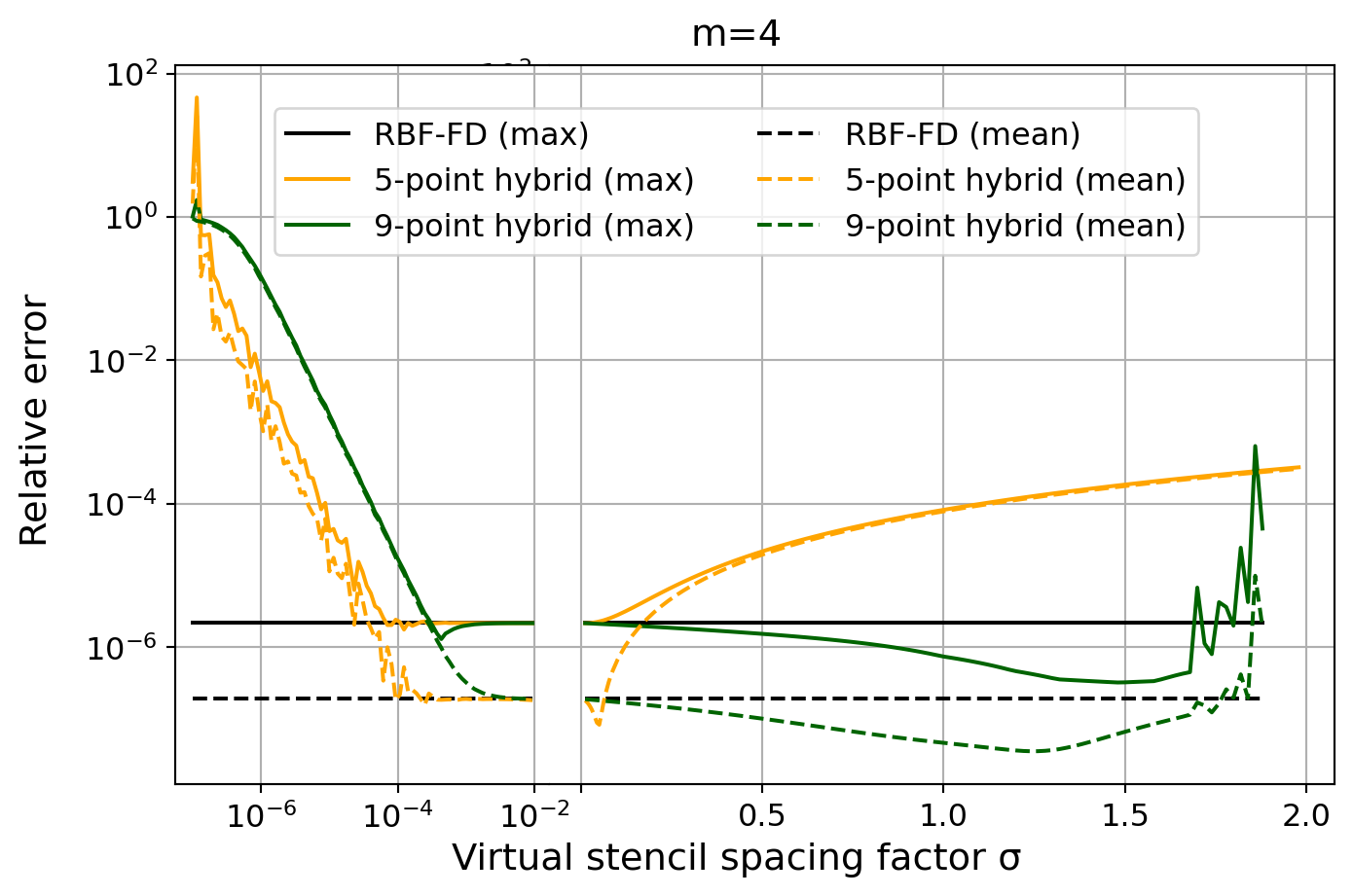}
    \caption{Error comparison for 4th order monomial augmentation and the five-point and nine-point stencil.}
    \label{fig:error-m4}
\end{figure}

Next, we analyze the augmentation of order $m=4$. Fig. \ref{fig:error-m4} again compares the errors for the five-point and nine-point stencils. While the nine-point stencil achieves greater accuracy than RBF-FD for a certain range of $\sigma$, the difference is not as big as in the $m=2$ case. The five-point stencil generally performs worse than pure RBF-FD, as is expected, since in this case, the error is dominated by a lower order finite difference stencil. However, surprisingly, for very small values of $\sigma$ that is no longer the case and the error behavior is similar to that observed for the second-order augmentation in Fig. \ref{fig:error-m2}, with the start of the increase in error at values of $\sigma$ about an order of magnitude larger. 

\begin{figure*}
  \centering
  \includegraphics[width=0.49\linewidth]{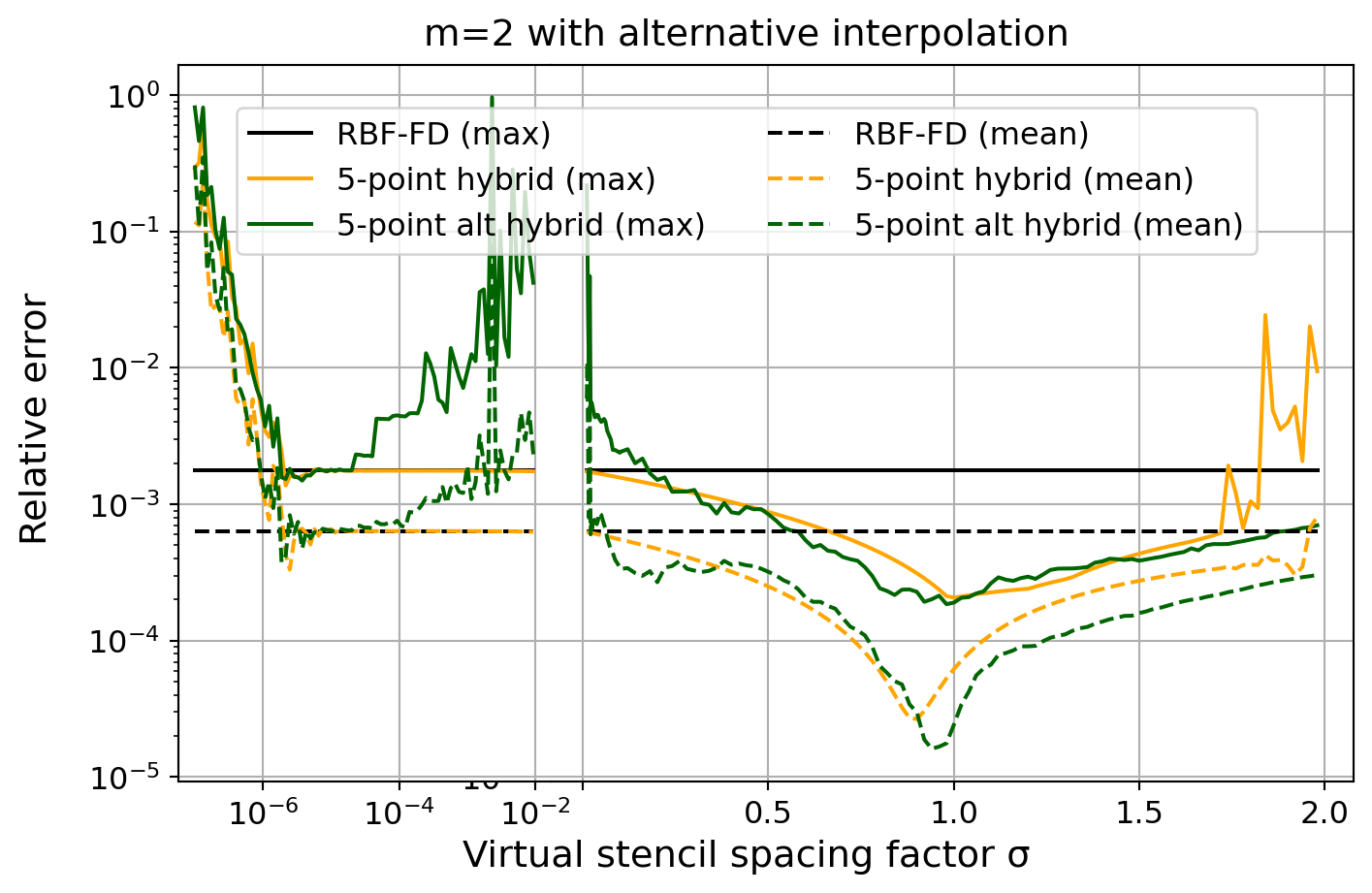}
  \includegraphics[width=0.49\linewidth]{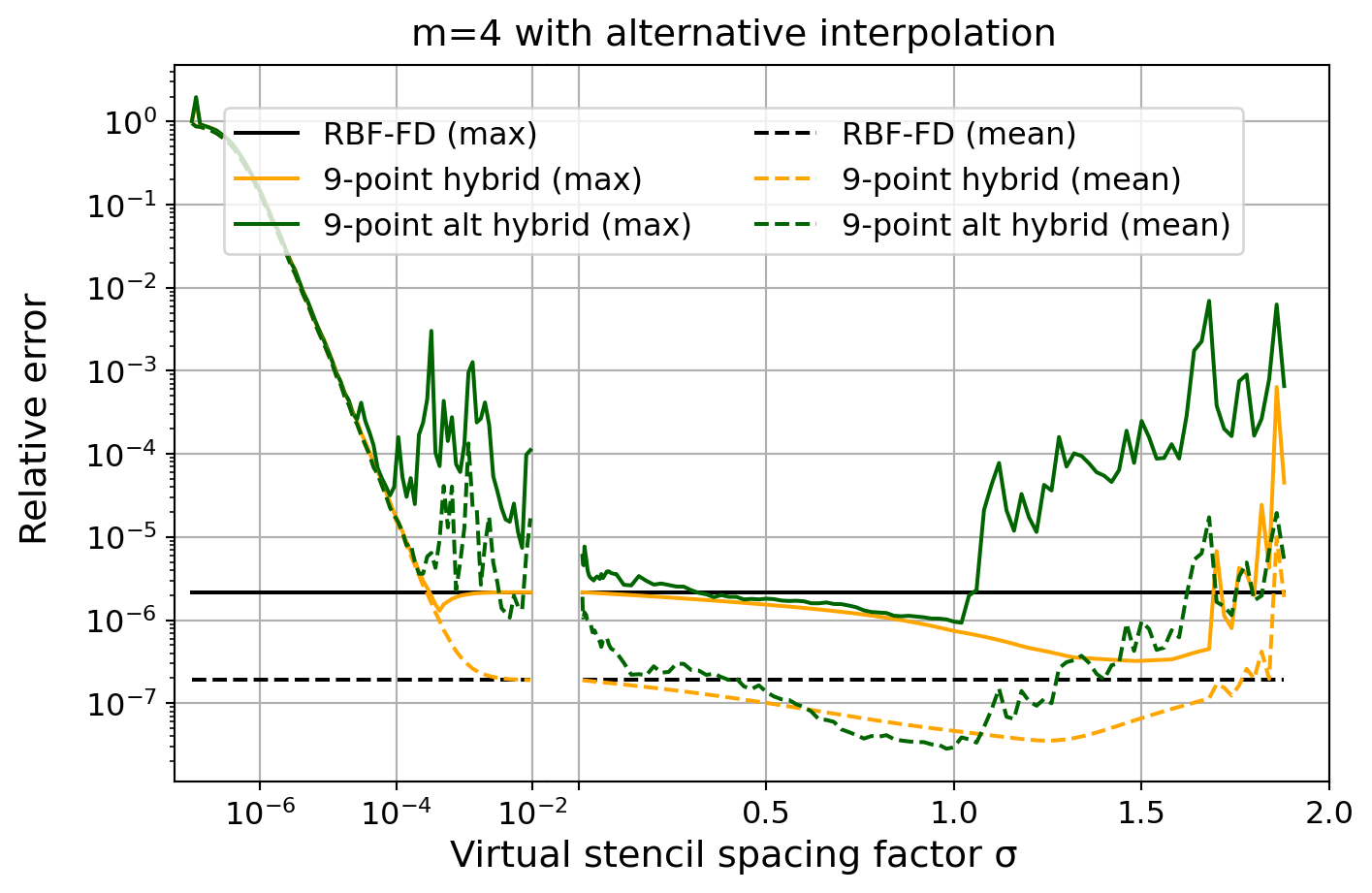}
  \caption{Error comparison for the alternative interpolation method.}
  \label{fig:error-alternative}
\end{figure*}

In the cases considered so far, the error starts to increase as we further increase $\sigma$ away from one, which, as mentioned, might be because virtual stencil nodes become too distant from the stencil center, resulting in an ineffective interpolation. This leads us to consider an alternative approach: instead of always using the same stencil of the center point $\boldsymbol{x}$ to interpolate to all the virtual stencil positions as in (\ref{eqn:hybrid-interp}), we can use different stencils for each of the virtual nodes $\boldsymbol{x} + \Delta_i$, namely, the $n$ closest nodes to the given virtual node, as opposed to the $n$ closest nodes to the center point $\boldsymbol{x}$. We investigated the effect of this approach for $m=2$ with the five-point stencil and for $m=4$ with the nine-point stencil. The result are presented on Fig. \ref{fig:error-alternative}. Surprisingly, the alternative approach does not result in any noticeable improvement while being much more computationally expensive, as we need to set up and solve a completely separate linear system for each node from the virtual stencil and can no longer benefit from saving the LU decomposition of the matrices.

We now compare the computation time of the hybrid method and RBF-FD. Each configuration was run 25 times on an Intel® Core™ i5-8250U processor and we took the median runtime. While we have verified that both methods have the same cubic asymptotic growth in the first phase with respect to the number of support nodes $n$, this overlooks the actual runtime differences at lower, relevant values of $n$, where lower-order terms still have a noticeable impact, as shown in Table \ref{tab:computation-time-phase1}. In the second phase, the resulting large sparse linear system is solved. In practice, this is done using iterative methods. The number of iterations required by an iterative solver is difficult to predict, as it depends on the structure and conditioning of the matrix. For our analysis, we used the biconjugate gradient stabilized iterative method (BiCGSTAB) preconditioned by an incomplete LU factorization with thresholding (ILUT). We observed no significant differences in runtime (see Table \ref{tab:computation-time-phase2}), except for larger values of $\sigma$, where the error and iteration count increased substantially.

\begin{table}
    \centering
    \begin{tabular}{c|c|c|c}
        Configuration & RBF-FD & 5-point hybrid & 9-point hybrid \\ \hline
        $m=2, n=12$ & 41 & 57 & 79 \\ \hline
        $m=4, n=30$ & 218 & 268 & 331 \\ \hline
        $m=2, n=96$ & 1323 & 1473 & 1672 \\
    \end{tabular}
    \caption{Computation time in milliseconds for the first phase, $\sigma = 1$.}
    \label{tab:computation-time-phase1}
\end{table}

\begin{table}
    \centering
    \begin{tabular}{c|c|c|c}
        Configuration & RBF-FD & 5-point hybrid & 9-point hybrid \\ \hline
        $m=2, n=12$ & 130 & 130 & 130 \\ \hline
        $m=4, n=30$ & 612 & 607 & 614 \\ \hline
        $m=2, n=96$ & 4249 & 4210 & 4219 \\
    \end{tabular}
    \caption{Computation time in milliseconds for the second phase, $\sigma = 1$.}
    \label{tab:computation-time-phase2}
\end{table}

\begin{figure}
  \centering
  \includegraphics[width=\linewidth]{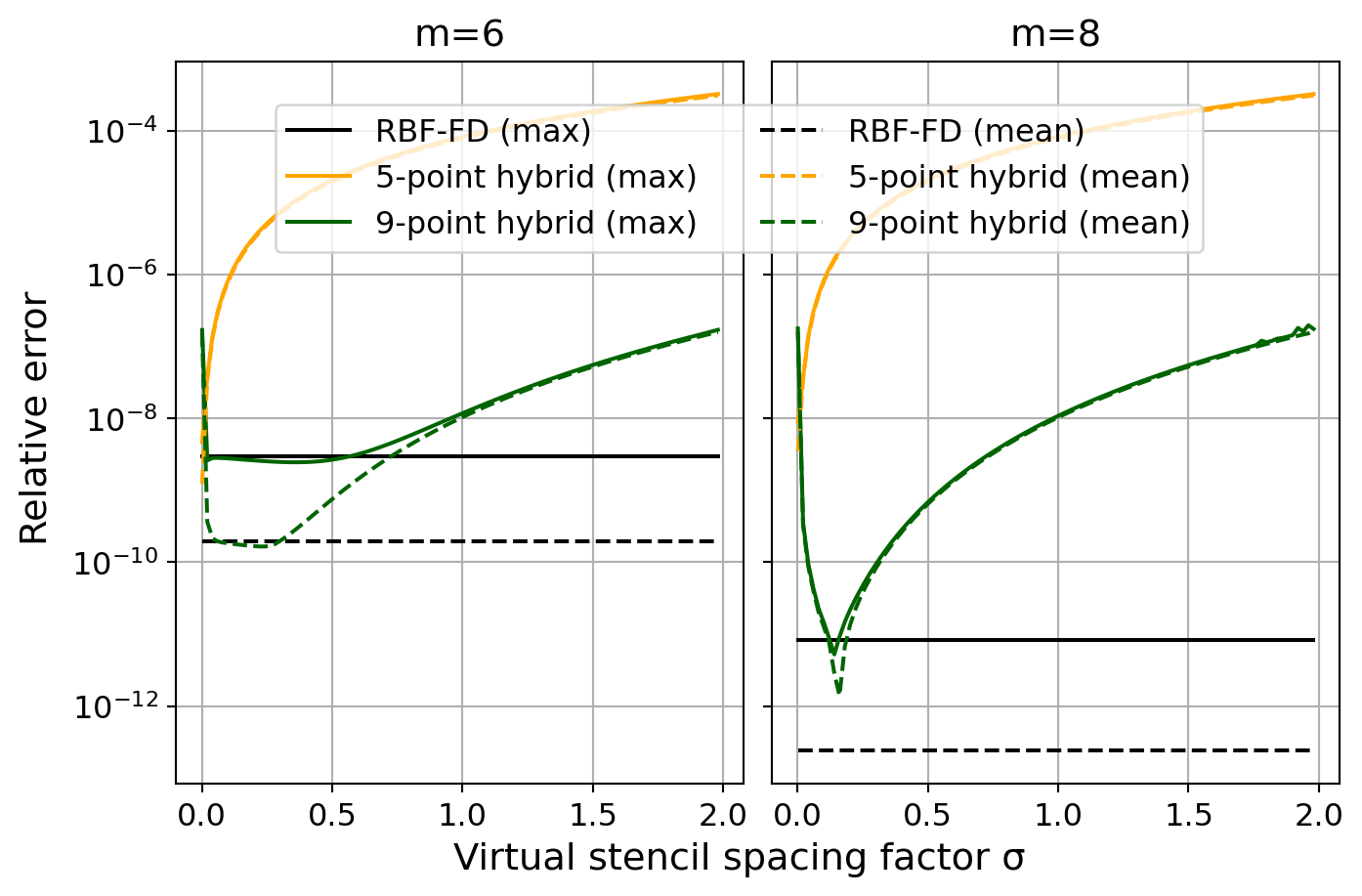}
  \caption{Error comparison for 6th and 8th order monomial augmentation.}
  \label{fig:error-m6-m8}
\end{figure}

Lastly, we investigate how the hybrid method performs for even higher orders of monomial augmentation. Fig. \ref{fig:error-m6-m8} illustrates the cases of $m=6$ and $m=8$. Interestingly, we can fine tune the value of $\sigma$ so that the hybrid approach performs just as well as RBF-FD for $m=6$, even though the nine-point stencil only has a polynomial reproduction degree of $5$. For $m=8$ the RBF-FD method outperforms the hybrid approach for all tested values of $\sigma$ and the two finite difference stencils.

\section{Conclusion}

In this study, we compared two RBF-based methods for approximating differential operators - the usual RBF-FD and a hybrid numerical method combining RBF interpolation with finite difference stencils. For both, we opted for a popular choice of polyharmonic spline RBFs with monomial augmentation.

We have compared the accuracy of the two approaches on a two-dimensional Poisson problem. 
Our results indicate that the hybrid method can achieve improved accuracy over RBF-FD, particularly when the virtual stencil step size aligns well with the discretisation distance.
However, this improvement comes at the cost of increased computational complexity due to the additional interpolation steps.
Furthermore, we noticed that by fine-tuning the virtual stencil size, the hybrid approach can benefit from a greatly increased accuracy, in some cases even matching a pure RBF-FD approach of a higher polynomial reproduction degree.

An important remark to be made is that our analyses have only been performed on a specific and relatively simple Poisson problem. As part of our ongoing research we plan to focus on more complex problems, employing some problem-specific finite difference stencils in the hybrid approach, specifically ones arising from a staggered grid finite difference approach to hyperbolic systems of PDEs.

\section*{Acknowledgement}
The authors would like to acknowledge the financial support from the Slovenian Research and Innovation Agency (ARIS) in the framework of the research core funding No. P2-0095, The Young Researcher program PR-12347 and research project No. J2-3048.

\bibliographystyle{IEEEtran}
\bibliography{paper}

\end{document}